\documentclass{gtart_s}
\usepackage{amsmath,amssymb,mathrsfs}

\def\abs#1{\left| #1 \right|}
\def\Z{\mathbb Z}
\let\wtilde\widetilde
\let\what\widehat
 
\def\gen#1{\langle #1\rangle}

\theoremstyle{plain}
\newtheorem{theorem}{Theorem}

\newtheorem{lemma}{Lemma}
\newtheorem*{serre}{Serre's Theorem}

\theoremstyle{definition}
\newtheorem*{remark}{Remark}
\newtheorem*{remarks}{Remarks}

\begin{document}

\title{The adjunction problem and a theorem of Serre} 
\author{Max Forester\\Colin Rourke}

\address{Math Department, University of Oklahoma, Norman OK
73019, USA\\\smallskip\\Maths Institute, University of Warwick, 
Coventry, CV4 7AL, UK} 
\gtemail{\mailto{forester@math.ou.edu}{\rm\qua 
and\qua}\mailto{cpr@maths.warwick.ac.uk}}

\begin{abstract}
In this note we prove injectivity and relative asphericity for
``layered'' systems of equations over torsion-free groups, when the
exponent matrix is invertible over $\Z$.  We also give elementary
geometric proofs of results due to Bogley--Pride and Serre that are
used in the proof of the main theorem.
\end{abstract}

\makeshorttitle

A long-standing problem in group theory is the {\em
adjunction problem} of deciding when a given group injects into the
group obtained by adjoining new generators and relators
\cite{Neu}. This note solves the adjunction problem over torsion-free
groups in the special case in which new generators and relators are
added in pairs and such that the exponent matrix is invertible.  We
prove that in this case the group does inject.  The case of one such
pair was proved by Klyachko \cite{Kl}.  The extension uses our
previous paper \cite{FoR} and a result of Bogley and Pride
\cite{BoPr}, which in turn follows from an old theorem of Serre, whose
proof depends on Tate cohomology; see Huebschmann \cite{Hu}.  We give
new and elementary proofs of both the Bogley--Pride and Serre results.

Let $(L,K)$ be a relative $2$--complex (a CW--pair such that $L - K$ is
at most $2$--dimensional). We say that $(L,K)$ is \emph{relatively
aspherical} if the map $$\pi_2(K \cup L^{(1)},K) \to \pi_2(L,K)$$ is
surjective. As shown in \cite[3.1--3.3]{FoR}, this occurs if and only if
conclusions (a) and (b) of Theorem \ref{mainthm} below hold.  This is the
natural topological notion of asphericity but it should be noted that it
differs from the combinatorial notion introduced in \cite{BoPr}. The
difference concerns the definition of irreducibility of diagrams
representing elements of $\pi_2(L,K)$; see \cite{FoR}. 

The fundamental group of $L$ is obtained from $G = \pi_1(K)$ by adding
generators $\{t_i\}$ and relators $\{r_j\}$ corresponding to the
$1$--cells and $2$--cells respectively of $L-K$. The relators $r_j \in G
\ast \langle t_1, \ldots, t_n\rangle$ can then be viewed as a system of
equations in the variables $\{t_i\}$ with coefficients in $G$. It is well
known (see Howie \cite{Ho} for example) that the map $\pi_1(K) \to \pi_1(L)$ is
injective if and only if the system has a solution in an overgroup of
$G$. 

The \emph{exponent matrix} of the system (or of the pair $(L,K)$) has
entries $m_{ij}$ equal to the exponent sum of $t_i$ in the relator
$r_j$. In topological terms it is the $2$--dimensional boundary map in
the relative cellular chain complex of $(L,K)$. 

A long-standing conjecture \cite{Ho} states that for any relative
$2$--complex $(L,K)$, if the exponent matrix is nonsingular, then
$\pi_1(K) \to \pi_1(L)$ is injective. If we assume further that
$\pi_1(K)$ is torsion-free and the exponent matrix is invertible over
$\Z$ then we conjecture that $(L,K)$ is also relatively
aspherical. (As shown in \cite{FoR} this conclusion can fail if either
of the additional hypotheses is omitted.)  Our main result proves this
in a special case:

\begin{theorem}\label{mainthm}
Let $(L,K)$ be a layered relative $2$--complex with $\pi_1(K)$
torsion-free. If the exponent matrix is invertible over $\Z$ then
\begin{enumerate} 
\item[\textup{(a)}] $\pi_1 (K) \to \pi_1 (L)$ is injective, and
\item[\textup{(b)}] the inclusion-induced map $\Z\pi_1 (L)
\otimes_{\Z\pi_1(K)} \pi_2 (K) \to \pi_2 (L)$ is an isomorphism. 
\end{enumerate}
\end{theorem}

Here, $(L,K)$ is \emph{layered} if $L-K$ has equal numbers of $1$-- and
$2$--cells and $L$ is formed from $K$ by alternately adding 
$1$--cells and $2$--cells. In terms of the associated relative
presentation it means that the generators and relators can be added
alternately. 

A special case of Theorem \ref{mainthm} was proved in \cite{FoR}: the
theorem was proved when $L - K$ consists of one $1$--cell and one
$2$--cell.  In this note we observe that the special case can be
applied inductively, with the aid of the Bogley--Pride--Serre result,
stated and proved as Theorem \ref{bpsthm} below. It is worth 
stressing that part (b) of Theorem \ref{mainthm}, for the case of one
new generator and one new relator, is a non-trivial extension of
Klyachko's theorem and it is the key to allowing the inductive argument
of this note to proceed. 

\begin{proof}
The layered hypothesis implies that there is a nested sequence of
subcomplexes $K = K_0 \subset K_1 \subset \cdots \subset K_n = L$ where
$K_{i+1} - K_i$ has one $1$--cell and one $2$--cell for each $i$. Note
that the exponent matrix for $(L,K)$ is triangular with diagonal entries
equal to $\pm 1$, and these diagonal entries represent the $1\times 1$
exponent matrices for the pairs $(K_{i+1},K_i)$. In particular each pair 
$(K_{i+1},K_i)$ is \emph{amenable} in the sense of Fenn and Rourke
\cite{FR}.  Then by the main theorem of \cite{FoR}, the pair
$(K_{i+1},K_i)$ is relatively aspherical provided $\pi_1(K_i)$ is
torsion-free.

We are given that $\pi_1(K_0)$ is torsion-free, and so $(K_1,K_0)$ is
relatively aspherical. By Theorem \ref{bpsthm} below $\pi_1(K_1)$ is then
torsion-free. Proceeding inductively, using \cite{FoR} and Theorem
\ref{bpsthm}, we find that every pair $(K_{i+1},K_i)$ is relatively
aspherical. 

It remains to verify that relative asphericity is transitive. Given $K
\subset L \subset M$ with $M - K$ at most $2$--dimensional, relative
asphericity of $(M,L)$ and $(L,K)$ implies 
\begin{equation*}
\begin{split}
\pi_2(M) \  &= \ \Z\pi_1(M)
\otimes_{\Z\pi_1(L)} (\Z\pi_1(L) \otimes_{\Z\pi_1(K)} \pi_2(K)) \\ &=
\ \Z\pi_1(M) \otimes_{\Z\pi_1(K)} \pi_2(K)
\end{split}
\end{equation*}
 so condition (b) holds for
$(M,K)$. Condition (a) for $(M,K)$ is clear. 
\end{proof}

\begin{remark}
The proof shows that the exponent matrix hypothesis can be relaxed to
allow layered relative $2$--complexes for which each pair
$(K_{i+1},K_i)$ is amenable, ie, that the relator given by the new
2--cell has an ``amenable $t$--shape'' in terms of the new generator;
see \cite{FR} or \cite{FoR}.  The result also solves the adjunction
problem for systems of generators and relators which can be
transformed, by a change of variables, into a layered amenable system.
\end{remark}

The next theorem is analogous to (and follows from) Theorem 1.4 of
\cite{BoPr}, though our version also follows directly from the theorem
of Serre given in \cite{Hu} and used by \cite{BoPr}.  Serre's argument
is algebraic.  We shall give a direct geometric proof.  At the end of
the paper we extend this proof to recover the full power of Serre's
theorem.

\begin{theorem}\label{bpsthm}
If $(L,K)$ is relatively aspherical then every finite subgroup of
$\pi_1(L)$ is contained in a unique conjugate of $\pi_1 (K)$. 
\end{theorem}

\proof By adding cells of dimension $\geqslant 3$ we can arrange that
all the homotopy groups of $K$ vanish in dimensions 2 and above.  This
does not change the fact that $(L,K)$ is relatively aspherical.  The
easiest way to see this is to use the diagram interpretation used in 
\cite{FoR}: relative asphericity means that there are no irreducible
diagrams over $\pi_1(K)$ using the cells of $L-K$. This only depends on
$\pi_1(K)$ and the form of the added relators and hence is unchanged by a
change in the higher homotopy groups of $K$.  After adding the new cells
$\pi_2(L)$ is trivial. 

Let $\wtilde{L}$ be the universal cover of $L$ and $\wtilde{K}$ the
preimage of $K$ in $\wtilde{L}$. Let $\what{L}$ be the $2$--complex
obtained from $\wtilde{L}$ by collapsing each connected component of
$\wtilde{K}$ to a vertex. Since each of these components is contractible,
the map $\wtilde{L} \to \what{L}$ is a homotopy equivalence, and so 
$\pi_1(\what{L})$ and $\pi_2 (\what{L})$ are trivial. Then since $\what{L}$
is $2$--dimensional, it is contractible. 

Note that the induced action of $\pi_1 (L)$ on $\what{L}$ is free away from
the $0$--skeleton, and the vertices have stabilisers equal to the
conjugates of $\pi_1 (K)$ in $\pi_1 (L)$. Hence it suffices to show that
every finite subgroup of $\pi_1 (L)$ fixes a unique vertex of $\what{L}$. 

This follows from the next two lemmas, the first of which well-known,
cf \cite{Brown}.

\begin{lemma}\label{cyclic-case} Suppose a group acts cellularly on 
an acyclic finite dimensional CW--complex $Q$ freely away from the
$0$--skeleton. Then every non-trivial element $g$ of finite order
fixes a unique vertex of $Q$.
\end{lemma} 

\proof Denote $Q/\gen g$ by $T$ and let $f\co Q\to T$ be natural
projection.

{\bf Step 1}\qua{\sl If $g$ has prime order $p$ then $g$ has
at least one fixed point.}

Suppose that $g$ has no fixed points.  We will inductively construct
$\Z/p$--cycles $c_i$ in $T$ in all dimensions.  Start with $c_0$ any
vertex and let $b_0= f^{-1}c_0$.  Then $b_0$ consists of $p$ points
and hence is zero in $H_0(Q,\Z/p)= \Z/p$.  So $b_0$ is the boundary of
a 1--chain $a_1$ and we define $c_1=f(a_1)$.  Now suppose that $c_i$
has been constructed.  Let $b_i= f^{-1}c_i$ which is a $\Z/p$--cycle
in $Q$ and which $p$--fold covers $c_i$.  Since $Q$ is acyclic, $b_i$
is the boundary of a $\Z/p$--chain $a_{i+1}$ say.  Then
$c_{i+1}=f(a_{i+1})$ is the next cycle.

We claim that all these cycles are non-zero in $\Z/p$--homology.  It
then follows that $T$ (and hence $Q$) is infinite dimensional,
contradicting the hypotheses.  Hence $g$ has a fixed point.

To see the cycles are all essential notice that the construction is
natural and maps to a similar construction in the universal
$\Z/p$--bundle.  We use Milnor's construction for the universal bundle, 
namely $E = {\rm lim}_i *_iP$ where $P$ is a $p$ point space, $*_i$
denotes the $i$--fold join $P*P*\cdots*P$ (that is, $i$ join
operations on $i+1$ copies of $P$), and the action is the join of the
cyclic action on $P$.  If we apply the construction to the $i$-th
stage $*_iP\to R=*_iP\,/\,\Z/p$ then the cycles $b_j$ are the subsets
$*_jP$ for $j\le i$.  So in this case $b_i$ is the top (fundamental)
cycle in $*_iP$ and is therefore non-zero.  But if any of the cycles
$c_j$, $j\le i$ is zero, so are all subsequent ones and then $b_i$
would be zero.  It follows that $c_j$, $j\le i$ are non-zero in $R$
and hence, in the limit all the $c_i$ are non-zero.

\medskip
{\bf Step 2}\qua{\sl If $g$ has prime order $p$
then $g$ has at most one fixed point.} 

Suppose $g$ fixes at least two points (which must be vertices).
Choose two $x,y$ say.  Let $c_1$ be an arc in $T$ from $f(x)$ to $f(y)$
and let $b_1=f^{-1}(c_1)$ (a $\Z/p$--cycle) then $b_1$ bounds a chain
$a_2$ in $Q$.  Let $c_2=f(a_2)$ a $\Z/p$--cycle in $T$.  The
construction now proceeds as in step 1.  

We claim that as before all these cycles in $T$ are non-zero in
$\Z/p$--homology.  It then follows that $T$ (and hence $Q$) is
infinite dimensional, contradicting the hypotheses.  Hence $g$ has at
most one fixed point.

To see this consider the universal bundle $E$ as before.  We map $Q\to
\Sigma(E)$ (the suspension of $E$) by mapping $x$ and $y$ to the two
suspension points and any other fixed points to either suspension
point.  Now for each fixed point $a$ let $A$ be its link in $Q$.   By
universality choose an equivariant map $A\to E$.  Then a neighbourhood
of $a$ is mapped conically.  Finally the map is extended to map the rest
of $Q$ to $E$ by universality.  Then since the construction is again
natural it maps to a similar construction in $\Sigma(E)$.  But here we
are constructing the suspensions of the classes constructed in step 1
which are all non-zero.

\medskip
To finish the proof of the lemma, suppose $g$ has order $n$ and fixes
no point of $Q$. If $g^k$ fixes a point $x$ for some $k>1$, then $g^k$
also fixes $gx$ (which is not $x$), and then a suitable power of $g$
contradicts step 2.  Hence $\gen g$ acts freely on $Q$ and then a
power of $g$ contradicts step 1.  Similarly if $g$ has at least two
fixed points, then a power of $g$ contradicts step 2.  \endproof

\begin{lemma}\label{gen-case}
If a finite group $G$ acts on a set $X$ in such a way that each
non-trivial element fixes a unique point, then $G$ has a global fixed
point.\end{lemma}

\proof Let $x_g$ denote the unique fixed point of $g$.  Note that
$h(x_g)= x_{hgh^{-1}}$ and hence $G$ acts on $\{x_g \mid g\in G-
\{1\}\}$.  So without loss we may assume that $X = \{x_g \mid g\in G-
\{1\}\}$.

Denote the stabilizer of $x \in X$ by $G_x$.  Choose any $x \in X$ and
let $\mathscr{O}\subset X$ be the orbit containing $x$, and let
$n=\abs{\mathscr{O}}$.  By the orbit stabiliser theorem
$\abs{G}=n\abs{G_x}$.  Notice that if $y\in \mathscr{O}$ then $G_x$ and 
$G_y$ are conjugate and hence $\abs{G_x}=\abs{G_y}$.  Notice also that
the hypothesis of unique fixed points implies that if $x\ne y$ then
$G_x-\{1\}$ and $G_y-\{1\}$ are disjoint.

Now define $S = \{g\in G-\{1\} \mid x_g \in \mathscr{O}\}$.  Then $S$
is the union over $y\in\mathscr{O}$ of the disjoint sets $G_y-\{1\}$ and 
we have $\abs{S} = n(\abs{G_x} - 1)$, which implies that $\abs{S}
>\frac{1}{2}(\abs{G}-1)$. Since $\mathscr{O}$ was an arbitrary orbit
there is not room for another such set $S$ and we must have $S =
G-\{1\}$.  Thus $\abs{G}-1=n(\abs{G_x} - 1)$ which implies $n=1$ and
$G_x=G$. This completes the proof of Lemma~\ref{gen-case}, and of
Theorem~\ref{bpsthm}.  \endproof

\begin{remark}
If one is interested only in the case where $\pi_1(K)$ is
torsion-free, then only step 1 of Lemma \ref{cyclic-case} is needed.
For if $\pi_1(L)$ has torsion, then there is an element $g$ of prime
order, and it fixes a vertex as in step $1$.  But then $\pi_1(K)$ has
torsion, a contradiction.  Thus for the proof of Theorem \ref{mainthm}
the foregoing proof can be shortened to a few lines. 
\end{remark}

We now extend the proof just given to give a full proof of Serre's
Theorem.

First note that the proof of Lemma \ref{cyclic-case} did not use the
hypothesis that $Q$ is finite dimensional, but that each quotient
$Q/\gen g$ had finite $\Z/p$--homological dimension for each element
of prime order $p$ in $G$.  Nor was the hypothesis that $Q$ is acyclic
used fully.  Lemma \ref{gen-case} used nothing about $Q$.  Thus putting
the two proofs together we have the following. 

\begin{theorem}[Global Fixed Point Theorem]\label{gfp}Suppose that a 
finite group $F$ acts cellularly on a CW--complex $Q$ freely away from
the $0$--skeleton.  Suppose further that
\begin{enumerate} 
\item[\textup{(a)}]for each prime factor $p$ of $|F|$, $Q$ is
$\Z/p$--acyclic, and 
\item[\textup{(b)}]for each element $g$ of $F$ of prime order $p$ the
quotient $Q/\gen g$ has finite $\Z/p$--homological dimension.
\end{enumerate}
Then $F$ has a unique fixed point.\end{theorem}

We take the statement of Serre's Theorem from Huebschmann \cite{Hu}.

\begin{serre}
Let $G$ be a group and $\{G_i\}_{i\in I}$ a family of subgroups such
that for every $q\ge q_0$ the canonical map $H^q(G,M)\to \prod_i
H^q(G_i,M)$ is an isomorphism for every $G$--module $M$.  Then each
finite subgroup $F$ of $G$ is contained uniquely in a conjugate of
one of the $G_i$ (and does not meet any other such conjugate).
\end{serre}

\proof Let $K$ be the disjoint union of the $K(G_i,1)$ for $i\in I$
and form the open wedge $K^+$ (ie add an arc to an external basepoint
for each component) and then construct $L$, a $K(G,1)$, by attaching
cells to $K^+$.

Let be $\wtilde L$ universal cover of $L$ and $\wtilde K$ the inverse
image of $K$ in $\wtilde L$ (which comprises a number of copies of
universal covers of the $K(G_i,1)$'s).  Then form $\what L$ by
squeezing each component of $\wtilde K$ to a point.  Then since we are
squeezing contractible subcomplexes $\what L$ is contractible and $G$
acts freely off the $0$--skeleton.  Further the stabilisers of the
vertices are the conjugates of the subgroups $G_i$, so we have to
prove that each finite subgroup $F$ of $G$ fixes a unique vertex.

To do this we use Theorem \ref{gfp} with $Q=\what L$.  The space $\what
L$ is contractible, so we have hypothesis (a).  We have to check (b).

Now if $H$ is a subgroup of $G$ then $\what L/H$ is formed from a
cover of $L$ by squeezing components of the preimage of $K$.  But the
cohomology hypotheses lift to any cover (since they are given ``for
any $G$--module") so $\what L/H$ is formed by squeezing a subspace
which carries all but finitely many of the cohomology groups and hence
by excision it has finite (co)homological dimension.  Thus we have 
hypothesis (b) of Theorem \ref{gfp}. \endproof

\begin{remarks}(1)\qua The Global Fixed Point Theorem is much stronger
than needed to prove Serre's theorem.  In this application $Q$ was
contractible instead of just $\Z/p$--acyclic for certain $p$ and $Q/\gen
g$ had finite homological dimension with all coefficients.

(2)\qua We did not need $Q$ to be a CW--complex for the above proofs
to work but merely that the inclusion of the fixed points is a
cofibration.
\end{remarks}

\end{document}